\documentclass[12pt]{amsart}
\usepackage{amsmath}
\numberwithin{equation}{section}
\usepackage{amsthm}
\usepackage{amssymb}
\usepackage{amsfonts}
\usepackage{pinlabel}
\usepackage{graphicx}
\newcommand{\nc}[2]{ \newcommand{#1}{#2} }

\nc{\avint}{ {- \hspace{-3.5mm} \int} }  

\nc{\R}{\mathrm{I \! R}}  
\nc{\N}{\mathrm{ I \! N}}  
\newcommand{\closure}[1]{ \stackrel{\rule{.1 in}{.01 in}}{#1} }

\newcommand{\qclosure}[1]{ \stackrel{\rule{.4 in}{.01 in}}{#1} }
\newcommand{\qhclosure}[1]{ \stackrel{\rule{.45 in}{.01 in}}{#1} }

\newcommand{\chisub}[1]{ {\mathbf{\chi}}_{_{#1}} }

\newcommand{\refeqn}[1]{ (\!\!~\ref{eq:#1}) } 
\newcommand{\refthm}[1]{ \!\!~\ref{#1} }    

\nc{\Holder}{H\"{o}lder\ }

\nc{\ith}{ \ensuremath{\text{i}^{\text{th}}} }
\nc{\jth}{ \ensuremath{\text{j}^{\text{th}}} }
\nc{\kth}{ \ensuremath{\text{k}^{\text{th}}} }
\nc{\dst}{ \ensuremath{\text{1}^{\text{st}}_{\delta}} }
\nc{\dnd}{ \ensuremath{\text{2}^{\text{nd}}_{\delta}} }
\nc{\ost}{ \ensuremath{\text{1}^{\text{st}}} }
\nc{\tnd}{ \ensuremath{\text{2}^{\text{nd}}} }
\nc{\curl}{ \nabla \times }
\nc{\Div}{ \nabla \cdot }
\nc{\DC}{K}

\newcommand{\BVPc}[4]{  
  \begin{equation}
        \begin{array}{rl}
           #1 & \ \text{in}
               \ \ #4 \vspace{.05in} \\
           #2 & \ \text{on} \ \ \partial #4 \;,
        \end{array}
  \label{eq:#3}
  \end{equation}    }

\nc{\Ppl}{ \mathcal{M}^{+} }  \nc{\Pmn}{ \mathcal{M}^{-} }

\nc{\smiley}{ $\stackrel{\because}{\smile} \;$ }

\newcommand{\dist}{\mathrm{dist}}

\renewcommand{\eqref}[1]{(\ref{#1})}
\newcommand{\eat}[1]{}

\newcommand{\bsmall}{\begin{array}[c]{c}}
\newcommand{\esmall}{\end{array}}

\newcommand{\inj}{\mathrm{inj}}

\newcommand{\Ric}{\mathrm{Ric}}

\theoremstyle{plain}
\newtheorem{theo}{Theorem}[section]
\newtheorem{lem}[theo]{Lemma}
\newtheorem{prop}[theo]{Proposition}
\newtheorem{coro}[theo]{Corollary}

\theoremstyle{definition}
\newtheorem{defi}[theo]{Definition}

\newtheorem{Assumptions}[theo]{Assumptions}

\newtheorem{rmk}[theo]{Remark}

\numberwithin{equation}{section}

\def\>{>_{\sigma}}

\title[Mean Value Theorems for Riemannian Manifolds]{Mean Value Theorems for Riemannian Manifolds via the Obstacle Problem}
\author[Benson]{Brian Benson}
\address{Kansas State University, Mathematics Department, 
138 Cardwell Hall, 
Manhattan, KS 66506}
\email{babenson@ksu.edu}
\author[Blank]{Ivan Blank}
\address{Kansas State University, Mathematics Department, 
138 Cardwell Hall, 
Manhattan, KS 66506}
 \email{blanki@math.ksu.edu}
\author[LeCrone]{Jeremy LeCrone}
 \address{University of Richmond, Department of Math \& Computer Science,
212 Jepson Hall,
28 Westhampton Way,
University of Richmond, VA 23173}
 \email{jlecrone@richmond.edu}

\begin{document}
\baselineskip 18pt

\begin{abstract}
We develop some of the basic theory for the obstacle problem on Riemannian Manifolds,
and we use it to establish a mean value theorem.  Our mean value theorem works for
a very wide class of Riemannian manifolds and has no weights at all within the integral.
\end{abstract}
\maketitle

\setcounter{section}{0}

\section{Introduction}   \label{Intro}

The mean value theorem (MVT) is a fundamental tool in the analysis of harmonic 
functions and elliptic PDEs. Recalling the elementary setting of harmonic functions
in Euclidean space, the MVT states that the value of a harmonic function at a point
can be exactly recovered by taking the average value of the function over any sphere, 
or ball, centered at the point of interest (i.e. the {\it mean value property}
holds for all harmonic functions). More generally, if the function is sub-- 
(or super--) harmonic, the MVT states that the integral average of the function over 
a sphere gives an under-- (or over--) estimate for the value of the function at the center;
see for instance \cite[Theorem 2.1]{GT}. In this setting, where the MVT becomes an
inequality, one considers whether or not 
integral averages are suitably monotone (as a function of sphere radius) and
whether the value of the function can be recovered from a limiting process.

The importance of the MVT for the theory of harmonic functions and elliptic PDEs
is readily apparent, as a fundamental tool in proving both the weak and
strong maximum principle, the Harnack inequality, and a priori estimates.
Thus, as research on harmonic functions turned from Euclidean spaces to 
Riemannian manifolds, much effort has gone into developing versions of the 
MVT in this more general setting (c.f. \cite{Fr,GyWl,N,STY}). 
Reviewing the variety of statements for the MVT on manifolds, it becomes clear that 
properties of the manifold itself (regularity, curvature, topology, etc.) directly
affect what type of results one can hope to derive.

From early work of Friedman \cite{Fr} and Willmore \cite{W}, the direct
translation of the mean value property to a Riemannian manifold (taking 
integral averages over geodesic balls, or spheres) requires restriction to
so--called {\it harmonic} manifolds. In particular, Willmore proved that the manifold's
volume density function is spherically symmetric if and only if the classical 
MVT holds for geodesic balls \cite{W}.
Weakening the equality condition,
and looking at convergence properties of integral averages (as radius goes to 
zero), further statements have been derived on Einstein manifolds, and manifolds
with specific curvature restrictions. These attempts at widening the class of 
manifolds on which MVT statements can be proved has been accompanied by 
further modification and restrictions on the inequality itself. These 
modifications include introduction of weights (both inside and outside the 
integral), mean value inequalities which do not become equalities in the limit,
further specific conditions on geometric and topological properties of manifolds,
and terms involving Green's functions and their derivatives. 
One recent version of MVT due to Lei Ni states that harmonic functions satisfy 
\begin{equation}
    u(x) = \frac{1}{r^n} \int_{\Omega_r} |\nabla \log G|^2 u \; d\mu
\label{eq:NiMVT}
\end{equation}
where we see derivatives of the Green's function appearing, and domains of
integration, $\Omega_r$, are upper level sets for the Green's function 
\cite[Theorem 2.3]{N}.

Addressing some of the drawbacks in the previous MVT results, we present in
this paper a version of MVT which contains no weights associated with integration,
works for a broad class of Riemannian manifolds, and becomes a perfect
equality in the case of harmonic functions. The domain of integration for
our statement is no longer geodesic balls, but rather sets, $D_{x_0}(r)$, 
related to solutions to an obstacle problem where the solution is constrained
to be below the Green's function.
This connection between MVT and the obstacle problem begins with an 
important observation of Luis Caffarelli.

In his Fermi Lectures on the Obstacle Problem, Luis Caffarelli 
gave an elegant proof of the MVT that did not rely on many of the basic symmetry 
and smoothness properties of the Laplacian \cite{C}.  Indeed, he asserted a statement
of the MVT holding for general second order uniformly elliptic divergence form
operators which was subsequently proven in complete detail by the second author
and Hao \cite{BH}.  This proof described by Caffarelli and completed by Blank and Hao used the obstacle
problem to create a key test function.  In the case of Laplace's equation on $\R^n,$ the
creation of this function is trivial because of the aforementioned properties of the Laplacian
on Euclidean space.  In this paper, we wish to extend this proof to the Laplace-Beltrami operator
on Riemannian manifolds, where the construction of the crucial test function is no longer trivial.

In terms of reproducing a MVT result, it is clear that harmonic 
(and even subharmonic) functions on compact manifolds
without boundary are uninteresting.  We will be interested in harmonic 
(and subharmonic) functions defined on
all noncompact manifolds, and those functions which are defined on
compact manifolds with boundary which we assume are strict submanifolds 
of an ambient manifold.  Here, in our usage
of the word ``strict,'' we mean that the complement of the 
original manifold contains an open set.

Our main theorem can now be stated:
\begin{theo}[Mean Value Theorem on Riemannian manifolds]  \label{RiemMVT}
Given a point $x_0$ in a complete Riemannian manifold
$\mathcal{M}$ (possibly with boundary), there exists a maximal number $r_0 > 0$ (which is
finite if $\mathcal{M}$ is compact) and a family of open sets $\{ D_{x_0}(r) \}$ for $0 < r < r_0,$
such that
\begin{itemize}
   \item[(A)] $0 < r < s < r_0$ implies, $\qhclosure{D_{x_0}(r)} \; \subset D_{x_0}(s),$ and 
   \item[(B)] $\lim_{r \downarrow 0} \dist_{x_0}(\partial D_{x_0}(r)) = 0,$ and
   \item[(C)] if $u$ is a subsolution of the Laplace-Beltrami equation, then
$$u(x_0) = \lim_{r \downarrow 0} \frac{1}{|D_{x_0}(r)|} \int_{D_{x_0}(r)} u(x) \; dx \;,$$
and $0 < r < s < r_0$ implies
$$\frac{1}{|D_{x_0}(r)|} \int_{D_{x_0}(r)} u(x) \; dx \leq \frac{1}{|D_{x_0}(s)|} \int_{D_{x_0}(s)} u(x) \; dx \;.$$
\end{itemize}
Furthermore, if $\mathcal{M}$ is a compact manifold with boundary and $x_0 \in \mathcal{M},$ then
\begin{equation}
     \lim_{r \uparrow r_0} \left[ \inf_{x \in \partial D_{x_0}(r)} \dist(x, \partial \mathcal{M}) \right] = 0 \;.
\label{eq:gettobdry}
\end{equation}
Finally, if $r < r_0,$ then the set $D_{x_0}(r)$ is uniquely determined as the noncontact set of
any one of a family of obstacle problems.
\end{theo}
\noindent
Because of the lack of weights involved, and the perfect equality for 
solutions, one can hope to find nice relationships between properties 
of the sets $D_{x_0}(r)$ and properties of the manifold. The $D_{x_0}(r)$
sets are characterized as solutions of appropriate obstacle problems, 
and indeed, in order to prove our main theorem, 
we must first develop some of the basic theory for the 
obstacle problem on a Riemannian manifold.

\subsection*{Acknowledgements} The authors thank Dave Auckly and Bob Burckel for helpful discussions.


\section{Setting, Terminology, and Geometric Estimates}   \label{SetTermGeom}

We will use the following basic notation and assumptions throughout the paper:
$$
\begin{array}{lll}
\mathcal{M} & \ & \text{a smooth connected Riemannian n-manifold} \\
\mathcal{N} & \ & \mathcal{N} \subset \mathcal{M} \ \text{and} \ \mathcal{N} \ 
\text{satisfies assumptions that we detail below} \\
g & \ & \text{the metric for our ambient manifold} \ \mathcal{M} \\
T_p\mathcal{M} & \ & \text{the tangent plane to} \ \mathcal{M} \ \text{at} \ p \\
\text{Vol}(S) & \ & \text{the volume of the set} \ S \\
\chisub{D} & \ & \text{the characteristic function of the set} \ D \\
\closure{D} & \ & \text{the closure of the set} \ D \\
\partial D & \ & \text{the boundary of the set} \ D \\
\Omega(w) & \ & \{x:  w(x) > 0 \} \\
\Lambda(w) & \ & \{x: w(x) = 0 \} \\
FB(w) & \ & \partial \Omega(w) \cap \partial \Lambda(w) \\
\dist_p(x) & \ & \text{the distance from} \ p \ \text{to} \ x \ \text{in} \ \mathcal{M} \\
B_{p}(r) & \ & \{ x \in \mathcal{M} : \dist_p(x) < r \} \\
\eta_{\delta}(S) & \ & \text{the} \ \delta \ \text{neighborhood of the set} \ S \\
\inj_p(\mathcal{N}) & \ & \text{the injectivity radius of} \ \mathcal{N} \ \text{at} \ p \\
\Ric_p(v,w) & \ & \text{Ricci curvature of vectors }v,w \in T_p\mathcal{M} \\
H_p(x) & \ & \text{mean curvature of the geodesic sphere } \partial B_p(\dist_p(x)) \ \text{at} \ x \\
& & \text{(with respect to {\it inward} pointing normal vectors)}\\
D_p(r) & \ & \text{the Mean Value ball given in Theorem\refthm{RiemMVT}} \\
\Delta_g & \ & \text{the Laplace-Beltrami operator on} \ \mathcal{M} \\
\end{array}
$$
Note that in order to have a subscript ``$p$'' always denote a base point, we have
switched the locations of the center and radius from what is most customary for the notation for a ball,
including what was used in \cite{BH} for the mean value sets.

\begin{Assumptions}
Regarding some additional assumptions on $\mathcal{N}:$
\begin{enumerate}
  \item We assume $\mathcal{N}$ is a connected compact n-manifold 
		with smooth boundary and positive injectivity radius, call it $R_0.$
  \item We assume that $\mathcal{N}$ has the property that given any
		$f \in C^{\infty}(\mathcal{N})$ and $h \in C^{\infty}(\partial \mathcal{N}),$
		there exists a unique solution to:
  \begin{equation}
        \begin{array}{rl}
           \Delta_g u = f & \ \text{in}
               \ \ \mathcal{N} \vspace{.05in} \\
           u = h & \ \text{on} \ \ \partial \mathcal{N} \;.
        \end{array}
  \label{eq:DPsolvable}
  \end{equation}
\end{enumerate}
\end{Assumptions}

\begin{rmk}[Boundaries of Balls]  \label{BoB}
Note that if $\mathcal{M}$ is a manifold with boundary, and if $B_{x_0}(r)$
has a radius $r > 0$ which is sufficiently large to guarantee that a portion of
$\partial \mathcal{M}$ is within $B_{x_0}(r),$ then $\partial B_{x_0}(r)$
not only contains all $x \in \mathcal{M}$ with $\dist_{x_0}(x) = r,$
but it also contains all $x \in \partial \mathcal{M}$ with
$\dist_{x_0}(x) < r.$
\end{rmk}

In the remainder of this section, we collect fundamental geometric estimates
on manifolds and important results we will require for analysis on manifolds.
The first important result is the following geometric estimate which describes
the behavior of the Laplace--Beltrami operator acting on the square of the 
geodesic distance function (an important tool we use to help establish necessary
growth estimates for solutions to obstacle problems in the next section).
\begin{theo}[Key Estimate]  \label{KE}
If $r > 0$ is sufficiently small, then for all $x \in B_p(r)$
	\begin{equation}  \label{eq:RatherKeyEqn}
		|\Delta_g (\dist^2_p(x)) - 2n + \frac{1}{3} \Ric_p(v,v)  r^2 | 
			\leq C_{_{KE_1}} r^3
	\end{equation}
where $C_{_{KE_1}}$ is a finite positive constant depending on the manifold, 
and $v$ is the unit vector in the tangent plane whose direction corresponds 
to the geodesic flow from $p$ to $x.$
In particular, it follows that, in the same setting, there is a 
$C_{_{KE_2}} > 0$ such that
	\begin{equation}  \label{eq:RKEsimpvers}
		|\Delta_g (\dist^2_p(x)) - 2n| \leq C_{_{KE_2}} r^2 \;.
	\end{equation}
\end{theo}
\noindent
To prove Theorem\refthm{KE}\!\!, 
we use the following lemma from Li \cite[Theorem 4.1]{PL}:	

\begin{lem} \label{LapDist}
If $p \in \mathcal{M}$ and $x$ is not in the cut--locus of $p$, then
	\begin{equation}\label{eq:MCdistance}
		\Delta_g (\dist_p (x)) = H_p(x) \ .
	\end{equation}
\end{lem}

\noindent \textit{Proof of Theorem\refthm{KE}\!.} Applying the product rule, we have 
\begin{equation}\label{eq:ProdRule}
		\Delta_g (\dist^2_p(x))= 2\dist_p(x) \Delta_g (\dist_p(x)) 
			+2\underbrace{\|\nabla \dist_p(x) \|^2}_{=1} = 
				2\dist_p(x) \Delta_g (\dist_p(x)) + 2.	
	\end{equation}
Therefore, we focus on estimating $\Delta_g (\dist_p(x)),$ for which we
employ Lemma \ref{LapDist}.

We begin by working in geodesic normal coordinates, parametrizing geodesic balls
$B_p(r)$ via the exponential map (see Gray \cite{Gray} and Gray and Vanhecke 
\cite{GrayVan} for related computations). In particular, given $0 < \rho < r$ and 
$x \in \partial B_p(\rho)$ (not in the cut--locus of $p$),
it follows that $x = \exp_p(\rho v)$ for a unique 
$v \in S^{n-1} \subset T_p\mathcal{M}$. The expansion of the mean curvature of a 
geodesic ball is well-known, appearing for instance in \cite{PX}, and can be written as 
	\begin{equation}  \label{eq:MeanTaylor}
		H_p(x)=\frac{n-1}{\rho} - \frac{1}{3}\Ric_p (v,v)\rho+O(\rho^2).
	\end{equation}

Thus, combining \eqref{eq:MeanTaylor} with Lemma \ref{LapDist} and 
\eqref{eq:ProdRule}, we have
	\begin{equation}  \label{eq:RatherKeyEqn2}
		\Delta_g (\dist^2_p(x)) 
			=2n-\frac{1}{3}\Ric_p(v,v)\rho^2+O(\rho^3),
	\end{equation}
and Theorem 2.2 is thus proved. 
\hfill $\square$
\begin{rmk}
It now follows that $H_p(x)>0$ for all $x \in \partial B_{p}(\rho)$ and 
all sufficiently small $\rho > 0$.
\end{rmk}

Another important tool from global analysis is the following Harnack inequality 
for positive harmonic functions; the statement is given in Li 
\cite[Theorem 6.1]{PL}, where it is attributed to Yau:
\begin{theo}[Harnack]  \label{theo:Harnack}
	Let $\mathcal{M}$ be a complete Riemannian $n$-manifold and let $p \in \mathcal{M}$. 
	Suppose the Ricci curvature on $B_{p}(2s)$ is bounded from below 
	by $-(n-1)K$ for some $K \geq 0$. Then for $u$ a positive harmonic
	function on $B_{p}(2s)$, there exists a constant $C_{_H} = C_{_H}(p,n)>0$,
	depending only on $p$ and $n$, such that 
		\begin{equation}\label{eq:Harnack}
				\sup_{B_p(s)} u \leq \left( \inf_{B_p(s)} u\right )
						\exp \Big (C_{_H} \left( 1 + s\sqrt{K} \right ) \Big ).
		\end{equation}
\end{theo}
	
Returning to Equation \eqref{eq:DPsolvable}, and 
viewing this problem as one of the calculus of variations, we can find
$u$ by minimizing
\begin{equation}
    D_f(u) := \int_{\mathcal{N}} |\nabla u|^2 + 2uf
\label{eq:CalcVarVer}
\end{equation}
among functions
\begin{equation}
   u \in \mathcal{K}_{h} := \{ w \in W^{1,2}(\mathcal{N}) \; : \; w - h \in W^{1,2}_{0}(\mathcal{N}) \} \;.
\label{eq:KhDef}
\end{equation}
A minimizer for this variational problem exists as long as the functional is bounded from below,
though some effort will be required to ensure the minimizer is indeed a solution
to Equation \eqref{eq:DPsolvable}. We address the equivalence of these settings in the remainder of 
this section.

Assume that we are given $f \in C^{\infty}(\mathcal{N})$ and $h \in C^{\infty}(\partial \mathcal{N})$
which satisfy
\begin{equation}
   0 < \lambda \leq f \leq \mu < \infty \ \ \text{and} \ \ h \geq 0.
\label{eq:fhbds}
\end{equation}
Given $\mathcal{N}$ and $f$ as above,
we seek a nonnegative function $w$ which satisfies the following semilinear boundary value problem:
  \begin{equation}
        \begin{array}{rl}
           \Delta_g w = \chisub{ \{ w > 0 \} }f & \ \text{in}
               \ \ \mathcal{N} \vspace{.05in} \\
           w = h & \ \text{on} \ \ \partial \mathcal{N} \;.
        \end{array}
  \label{eq:OPsemiform}
  \end{equation}
As mentioned earlier, we solve this boundary value problem by seeking
a minimizer of the functional
\begin{equation}
    D_f(w) := \int_{\mathcal{N}} |\nabla w|^2 + 2wf
\label{eq:OPCV}
\end{equation}
among
\begin{equation}
     w \in \mathcal{P}_{h} := \mathcal{K}_{h} \cap \{ w : w \geq 0 \},
\label{eq:PhDef}
\end{equation}
which is a convex subset of $W^{1,2}(\mathcal{N}).$

\begin{rmk}
We note that our proof of the Mean Value Theorem in Section~\ref{MVT} only requires
that one consider obstacle problems with $f \equiv const.$ We treat more general functions $f$ 
herein because the general obstacle problem on Riemannian manifolds has independent interest.
\end{rmk}

With our assumptions on $\mathcal{M}$ and $\mathcal{N},$ and because we have $D_f$ bounded from
below on a convex subset of a Hilbert space (by assumption), the following theorem
is now obvious using standard arguments from the calculus of variations:
\begin{theo}[Minimizers]  \label{ExistUnique}
There exists a unique minimizer of $D_f$ within $\mathcal{P}_{h}.$
\end{theo}

\noindent
Now by following the exact same procedure found within the second section of \cite{BH}
we can show the following:
\begin{theo}[Existence with a gap]  \label{AlmostExist}
With our assumptions as above, there is a pair of functions $W$ and $F$ such that
\begin{equation}
            \begin{array}{rl}
           \Delta_g W = F & \ \text{in}
               \ \ \mathcal{N} \vspace{.05in} \\
           W = h & \ \text{on} \ \ \partial \mathcal{N} \;,
        \end{array}
\label{eq:FBgap}
\end{equation}
where the function $F$ is nonnegative and in addition satisfies:
\begin{equation}
    \begin{array}{rll}
         F(x) &\!\! = 0 \ \ \ &\text{for} \ x \in \{ W = 0 \}^{\text{o}} \rule[-.1 in]{0 in}{.3 in} \\
         F(x) &\!\! = f(x) \ \ \ &\text{for} \ x \in \{ W > 0 \}^{\text{o}} \rule[-.1 in]{0 in}{.3 in} \\
         F(x) &\!\! \leq \mu \ \ \ &\text{for} \ x \in \partial \{ W = 0 \} \cap \partial \{ W > 0 \} \rule[-.1 in]{0 in}{.3 in} \;,
    \end{array}
\label{eq:Fprops}
\end{equation}
where for any set $S \subset \mathcal{N},$ we use $S^{\text{o}}$ to denote its interior.  Thus $F$ agrees with
$\chisub{ \{ W > 0 \} }f$ everywhere except possibly the free boundary, which is the set
$\partial \{ W = 0 \} \cap \partial \{ W > 0 \}.$
\end{theo}
In order to show that there is no gap (i.e. \;\!\!that $F = f$ a.e.), we follow the proof within \cite{BH}
which requires that we have the basic optimal regularity and nondegeneracy statements, and we
turn to the statement and proof of those statements now.

\section{Basic Theory of the Obstacle Problem}   \label{ObProb}

As we have mentioned, we wish to show that the free boundary has measure zero and so
we will actually have
\begin{equation}
   \Delta_g W = \chisub{ \{ W > 0 \} }f \;.
\label{eq:Goal}
\end{equation}
Within \cite{BH} there is a key regularity property and a key nondegeneracy property which lead to this
statement.   In order to prove these properties, we will adapt the proofs found in \cite{B} and/or \cite{C},
and in particular, we won't need anything as delicate as what is found in \cite{BH} for the
proof of nondegeneracy.  On the other hand, we will make essential use of the key estimates 
from Theorem\refthm{KE}\!.

\subsection{Optimal Regularity}  \label{OpRegS}

\begin{theo}[Local Optimal Regularity]  \label{OpReg}
Fix $p_0 \in \mathcal{N} \cap \partial \{ W > 0 \}.$
There exists an $r_{_{OR}} > 0$ and a $C_{_{OR}} > 0$ with a finite bound depending only on
$\dist_{p_0}(\partial \mathcal{N}),$ bounds on the Ricci curvature on $\mathcal{N},$
and on the constant $C_{_{KE_1}}$ from Equation \eqref{eq:RatherKeyEqn}, such that 
if $B_{p_0}(2r) \subset \mathcal{N}$ and $r \leq r_{_{OR}},$
then for any $x \in B_{p_0}(r)$ we have
\begin{equation}
   W(x) \leq C_{_{OR}} \, \mu \; \dist_{p_0}(x)^2
\label{eq:OpReg}
\end{equation}
where $\mu$ is the bound on $f$ from above given in Equation \eqref{eq:fhbds}.
\end{theo}
\noindent
We will also need the following corollary which we obtain by adjusting and 
iterating most of the argument found in the proof of the local optimal regularity result.
\begin{coro}[Boundedness Result]  \label{Bddness}
Fix $p_0 \in \mathcal{N} \cap \partial \{ W > 0 \},$ and fix a compact $\Gamma \subset \subset \mathcal{N}$
which contains $p_0.$
Then there is a constant $\delta(\Gamma)$ such that for any $x \in \Gamma$ we have
\begin{equation}
   W(x) \leq \delta(\Gamma) \;.
\label{eq:bddgam}
\end{equation}
\end{coro}

\noindent
\textit{Proof of Theorem\refthm{OpReg}\!\!.}  Multiplying through by $1/\mu$,
we can assume, without loss of generality, that $\mu = 1.$ We fix
\begin{equation}
    r_{_{OR}} := \min \left\{ \inj_{p_0}(\mathcal{N}), \sqrt{\frac{n}{C_{_{KE_2}}}} \right\}.
\label{eq:r0def}
\end{equation}
(Note that $C_{_{KE_2}}$ depends only on bounds on the Ricci curvature and on $C_{_{KE_1}}.$)
Consider $B_{p_0}(r),$ with radius $r \leq r_{_{OR}}.$
Let $C(p_0,s,K) :=\exp\Big (C_{_{H}}(p_0,n)\left (1+s\sqrt{K}\right ) \Big )$ 
where $C_{_{H}}(\cdot, \cdot)$ is the constant from Theorem\refthm{theo:Harnack}\!.
After observing that $C$ is monotone increasing with respect to $s,$ we
define $C_1 := C(p_0, r_0, K).$
Within the neighborhood $B_{p_0}(r),$ 
for any nonnegative
solution of the Laplace-Beltrami Equation $\Delta_g u = 0,$ the Harnack inequality from
Equation\refeqn{Harnack}in Theorem\refthm{theo:Harnack}gives:
\begin{equation}
     \sup_{B_{p_0}(s)} u \leq C_1 \inf_{B_{p_0}(s)} u
\label{eq:goodHarn}
\end{equation}
as long as $2s \leq r \leq r_{_{OR}}.$ Now, for such an $s,$ we let $u$ be the solution of:  
	\BVPc{\Delta_g u = 0}{u = W}{harmpart}{B_{p_0}(2s)}
we let $v$ be the solution of:  
	\BVPc{\Delta_g v = \Delta_g W}{v = 0}{nonharmpart}{B_{p_0}(2s)}
and we observe that $W = u + v.$  Obviously $0 \leq \Delta_g v \leq 1.$ Using
Equation\refeqn{RatherKeyEqn}we see that the function
\begin{equation}
   D(p) := \frac{\dist_{p_0}(p)^2 - 4s^2}{n},
\label{eq:Dpdef}
\end{equation}
satisfies:
\begin{equation}
    \Delta_g D(p) \approx 2,
\label{eq:Dpest}
\end{equation}
and now using the fact that $$r_{_{OR}} \leq \sqrt{\frac{n}{C_{_{KE_2}}}}$$ is implied by
Equation\refeqn{r0def}we have that $\Delta_g D \geq 1 \geq \Delta_g v.$
At this point, by noting that $v(p) = D(p) \equiv 0$ on $\partial B_{p_0}(2s),$ and
using the weak maximum principle, we get:
\begin{equation}
   \frac{-4s^2}{n} \leq D(p) \leq v(p) \leq 0
\label{eq:vest}
\end{equation}
within $B_{p_0}(2s).$ 
Now, 
we use the last equation 
to get
\begin{equation}
    u(p_0) = W(p_0) - v(p_0) = -v(p_0) \leq 4s^2/n \;.
\label{eq:UseZero}
\end{equation}
Applying the Harnack inequality in Equation\refeqn{Harnack}gives us:
\begin{alignat*}{1}
     \sup_{B_{p_0}(s)} u &\leq C_1 \inf_{B_{p_0}(s)} u \\
      &\leq C_1 u(p_0) \\
      &\leq 4C_1 s^2/n \;.
\end{alignat*}

From Equation\refeqn{vest}\!\!, we have $\sup_{B_{p_0}(s)} |v| \leq 4s^2/n.$ 
It follows that
	\begin{equation} \label{eq:absw} \begin{split}
		W(x) \leq |u(x)| + |v(x)| & \le 4C_1 s^2/n +4s^2/n \\
			&=4\big (C_1+1\big ) s^2/n \\
			&= \; (C_{_{OR}}) s^2	
	\end{split}
	\end{equation}
holds for all $x \in B_{p_0}(s),$ where $C_{_{OR}} := 4\big (C_1+1\big )/n$ and $2s \le r.$
We conclude the proof by fixing any $x \in B_{p_0}(r/2)$ and taking $s = \dist_{p_0}(x).$
Then\refeqn{absw}implies
\[
	W(x) \le \; C_{_{OR}} (\dist_{p_0}(x))^2 \;.
\]
\hfill $\square$

\noindent {\em Proof of Corollary\refthm{Bddness}\!\!.} The idea of this proof is that we iterate the
basic estimate of the previous proof until we have exhausted $\Gamma$ and, because $\Gamma$ is compact,
we get to every point within $\Gamma$ after a finite number of iterations.  The only real change occurs because
we have to iterate the estimate around points which are not part of the free boundary.  Indeed, if we have
$p_1$ near enough to $p_0$ so that we can apply the local result, and then we attempt to redo everything
centered at $p_1,$ then Equation\refeqn{UseZero}becomes
\begin{equation}
    u(p_1) = W(p_1) - v(p_1) = C_{_{PB}} - v(p_1) \leq C_{_{PB}} + 4s^2/n \;,
\label{eq:UsePB}
\end{equation}
where the ``PB'' in $C_{_{PB}}$ is short for the ``previous bound.''  This new estimate propagates into the new version of
Equation\refeqn{absw}\!\!, but we still get a finite bound on the right hand side.  The constant $C_{_{PB}}$ does
indeed get worse with each step, but is always finite, and we only need a finite number of steps to exhaust $\Gamma.$
\hfill $\square$

\subsection{Nondegeneracy} \label{NonDegen}

\begin{theo}[Local Nondegeneracy]  \label{NonDeg}
Fix $p \in \mathcal{N} \subset \subset \mathcal{M}.$
There exists an $r_{_{ND}} > 0$ depending only on
bounds on the Ricci curvature on $\mathcal{N}$
and on the constant $C_{_{KE_1}}$ from Equation \eqref{eq:RatherKeyEqn},
and there exists a constant $C_{_{ND}}$ depending only on $n$ such that 
if $p \in \; \qclosure{\Omega(W)}, \; B_{p}(r) \subset \mathcal{N},$ and $r \leq r_{_{ND}},$
then
\begin{equation}
    \sup_{x \in B_{p}(r)} W(x) \geq W(p) + C_{_{ND}} \lambda \, r^2
\label{eq:NonDeg}
\end{equation}
where $\lambda$ is the bound on $f$ from below given in Equation \eqref{eq:fhbds}.
\end{theo}
\noindent
We will also need the following global nondegeneracy statement which is easily obtained
by iterating parts of the argument from the proof of the previous theorem.
\begin{theo}[Global Nondegeneracy]  \label{GNonDeg}
Fix $p \in \; \qclosure{\Omega(W)}$ and
$\qclosure{B_{p}(R)} \; \subset \subset \mathcal{N}.$ There exists a constant $C_{_{GN}} > 0$
depending only on bounds on the Ricci curvature on $\mathcal{N},$
and on bounds over $\mathcal{N}$ on the constant $C_{_{KE_1}}$ from Equation \eqref{eq:RatherKeyEqn},
such that for any $r_{_{ND}} \leq s \leq R$ we have:
\begin{equation}
    \sup_{x \in B_{p}(s)} W(x) \geq W(p) + C_{_{GN}} s \, \lambda 
\label{eq:GNonDeg}
\end{equation}
where again the definition of $\lambda$ is found in Equation \eqref{eq:fhbds}.
\end{theo}

\begin{rmk}[Independence from $p$]  \label{IndepP}
If we are on a compact manifold, then $r_{_{OR}}$ and $r_{_{ND}}$ can be taken to be independent of $p$.
\end{rmk}

\noindent
\textit{Proof of Theorem\refthm{NonDeg}\!.}
Multiplying everything through by $1/\lambda$ we can assume that $f \geq 1$ without
loss of generality.  Next, by continuity of $w$ it suffices to take $p_1 \in \mathcal{N}$
satisfying $w(p_1) > 0$ and then prove
\begin{equation}
      \sup_{p \in \partial B_{p_1}(r)} w(p) \geq \frac{r^2}{2n+1} + w(p_1) \;.
\label{eq:supinside4p1}
\end{equation}
To this end, we define
\begin{equation}
     v(p):= w(p) - \frac{\dist_{p_1}(p)^2}{2n+1}.
\label{eq:KeyNDfct}
\end{equation}
By Equation\refeqn{RatherKeyEqn}if we take $r \leq r_{_{ND}}$ sufficiently small,
then for all $p \in B_{p_1}(r),$ we have:
\begin{equation}
    2n -1 \leq \Delta_g \dist_{p_1}(p)^2 \leq 2n + 1 \;.
\label{eq:dis2est}
\end{equation}
As a result, in $\Omega_r:=\{w >0\} \cap B_{p_1}(r)$, we have that 
	\begin{equation*}
		\Delta_g v(p) = 1 - \frac{\Delta_g \dist_{p_1}(p)^2}{2n+1} \geq 0
	\end{equation*}
and $v(p_1) = w(p_1).$
Therefore, since $v$ is subharmonic in $B_{p_1}(r)$, the weak maximum principle implies that 
	\begin{equation*}
		w(p_1) = v(p_1) \leq \sup_{p \in \partial B_{p_1}(r)} v(p)
		= \sup_{p \in \partial B_{p_1}(r)} w(p) - \frac{r^2}{2n+1} \;.
	\end{equation*}
Hence, 
	\begin{equation*}
		\sup_{p \in \partial B_{p_1}(r)} w(p) \geq \frac{r^2}{2n+1} + w(p_1) \;,
	\end{equation*}
from which the claim follows and we can take $$C_{_{ND}} := \frac{1}{2n + 1} \;.$$
\hfill $\square$

\noindent
\textit{Proof of Theorem\refthm{GNonDeg}\!.}  As in the previous proof, we can assume $f \geq 1.$
Now, we need to bound the function
\begin{equation}
     \Theta(s) := \sup_{x \in B_{p}(s)} ( W(x) - W(p) )
\label{eq:stufftoboundfrombelow}
\end{equation}
from below by a linear function of $s$ as long as $s \geq r_0.$  We already have
\begin{equation}
     \sup_{x \in B_{p}(r_{_{ND}})} ( W(x) - W(p) ) \geq C_{_{ND}} r_{_{ND}}^2 =: \gamma > 0
\label{eq:stuffwe have already}
\end{equation}
from local nondegeneracy.  We claim that
\begin{equation}
    \Theta(s) \geq \frac{\gamma}{2r_{_{ND}}} s
\label{eq:CartBlw}
\end{equation}
for all $s \geq r_{_{ND}}.$  Of course, for $r_{_{ND}} \leq s \leq 2r_{_{ND}}$ there is nothing to prove.

Now for $s \geq 2r_{_{ND}},$ we let $n=\left \lfloor \frac{s}{r_{_{ND}}}\right \rfloor,$ where
$\lfloor q \rfloor$ is the greatest integer less than or equal to $q.$  We set $p_0 := p$
and define $p_j$ iteratively for $j = 1 \ldots n$ by the equation:
\begin{equation}
     W(p_j) - W(p_{j - 1}) = \sup_{x \in B_{p_{j - 1}}(r_{_{ND}})} ( W(x) - W(p_{j-1}) ) \;,
\label{eq:itpjdef}
\end{equation}
and observe that 
by local nondegeneracy, this implies that
$W(p_j) - W(p_{j - 1}) \geq \gamma,$ and further, $p_j \in B_{p}(j r_{_{ND}}).$
Thus,
\begin{equation}
     \Theta(s) \geq \Theta(nr_{_{ND}}) \geq n\gamma
                     = \left \lfloor \frac{s}{r_{_{ND}}}\right \rfloor \gamma
                     \geq \frac{s}{2r_{_{ND}}} \gamma \;,
\label{eq:telescTheta}
\end{equation}
so we can let $$C_{_{GN}} := \frac{\gamma}{2r_{_{ND}}} \;.$$
\hfill $\square$

With the proof of optimal regularity and nondegeneracy behind us, we can use the exact same
proof as found in Lemma 5.1 of \cite{BT} and referenced within Corollary 3.10 of \cite{BH} to
show:
\begin{theo}[The Semilinear PDE Formulation]  \label{Upshot}
The set $\partial \{ W = 0 \} \cap \partial \{ W > 0 \}$ from Theorem\refthm{AlmostExist}is
strongly porous and therefore has a Hausdorff dimension strictly less than $n.$  As an immediate
consequence of this set having n-dimensional measure zero, we can say $W$ satisfies:
\begin{equation}
            \begin{array}{rl}
           \Delta_g W = \chisub{ \{ W > 0 \} }f & \ \text{in}
               \ \ N \vspace{.05in} \\
           W = h & \ \text{on} \ \ \partial N \;,
        \end{array}
\label{eq:nomoregap}
\end{equation}
and that solutions of this equation satisfy the regularity and nondegeneracy properties given
in Theorems\refthm{OpReg}and\refthm{NonDeg}above.
\end{theo}

\begin{rmk}
Now by following the same procedures as in the fourth section of Blank and Hao \cite{BH}, we can show that
the function $w$ from Theorem\refthm{ExistUnique}above is the same function found by the semilinear
PDE formulation in Equation\refeqn{nomoregap}(i.e. $w = W$).
\end{rmk}

\section{The Mean Value Theorem}   \label{MVT}

Now we have all of the tools that we need from the theory for the obstacle problem, and so
in this section we will split the proof of our main theorem into a number of smaller lemmas.
We start with the following.

\begin{lem}[Mean Value Theorem on Riemannian manifolds]  \label{RiemMVTagain}
Given a point $x_0$ in a complete Riemannian manifold,
$\mathcal{M},$ there exists a number $r_0 > 0$ and a family of open sets 
$\{ D_{x_0}(r) \}$ for $0 < r < r_0,$ such that
\begin{itemize}
   \item[(A)] $0 < r < s < r_0$ implies, $\qhclosure{D_{x_0}(r)} \; \subset D_{x_0}(s),$
   \item[(B)] $\lim_{r \downarrow 0} \dist(x_0, \partial D_{x_0}(r)) = 0,$ and
   \item[(C)] if $u$ is a subsolution of the Laplace-Beltrami equation, then
$$u(x_0) = \lim_{r \downarrow 0} \frac{1}{|D_{x_0}(r)|} \int_{D_{x_0}(r)} u(x) \; dx \;,$$
and $0 < r < s < r_0$ implies
$$\frac{1}{|D_{x_0}(r)|} \int_{D_{x_0}(r)} u(x) \; dx \leq \frac{1}{|D_{x_0}(s)|} \int_{D_{x_0}(s)} u(x) \; dx \;.$$
\end{itemize}
\end{lem}

\noindent
\textit{Proof of Lemma\refthm{RiemMVTagain}\!.}
We will first prove (C), and return to (A) and (B) afterward.
We will assume that $\mathcal{M}$ is not compact, but there is no real problem here as long as we
are confined to the case where $r_0$ is sufficiently small.
We fix $x_0 \in \mathcal{M},$ and we let $R > 1$ be very large and $r_0 > 0$ be very small.  
Let $B_{x_0}(R)$ be the open ball of radius $R$ centered at $x_0;$
or, in case $\partial B_{x_0}(R)$ is not smooth, we will assume that the ball
is approximated by a slightly larger set with smooth boundary.
Although these sets are not strictly geodesic balls, we will continue to refer to 
them as such for notational convenience. The reader will see in Lemma~\ref{UniqueMethod}
that this approximation has no meaningful effect on the mean value sets we ultimately construct
in the current proof.

Next, we let $G_R(x, x_{0})$ denote the Green's function on $B_{x_0}(R)$, with singularity
at $x_0$, and remind the reader that $G_R(x,x_0)$ automatically satisfies:
\begin{itemize}
    \item[I)] $G_R(x, x_0) > 0$ in $B_{x_0}(R) \setminus \{x_0\},$
    \item[II)] $G_R(x,x_0)$ is smooth in the set $B_{x_0}(R) \setminus \{x_0\},$
    \item[III)] $G_R(x,x_0) \equiv 0$ for all $x \in \partial B_{x_0}(R),$ and
    \item[IV)] $\lim_{x \rightarrow x_0} G_R(x,x_0) = + \infty.$
\end{itemize}
All of these properties are well-known, see for instance Aubin \cite[Chapter 4]{Aubin} 
and Li and Tam \cite[Section 1]{LT}.

Now we make definitions exactly like the ones above Theorem 4.1 of \cite{BH}.
We will fix $R$ and suppress this dependence for the rest of this argument until we potentially need to vary it
again.  Let
\begin{equation}
       J(w, B_{x_0}(R), r) := \int_{B_{x_0}(R)} |\nabla_g w|^2 - 2r^{-n} w \;,
\label{eq:fctldef}
\end{equation}
and let $w_r$ be the minimizer of $J$ among functions in $W^{1,2}_0(B_{x_0}(R))$ which are less than or equal to
the Green's function, $G_R(x,x_0).$  Note that the minimizer solves the obstacle problem:
\begin{equation}
      \begin{array}{rll}
             \Delta_g u &\!\!\!= -r^{-n}\chi_{_{ \{u < G_R \} }} \ \ \ & \text{in} \ B_{x_0}(R) \\
              \ \\
              u &\!\!\!= 0 \ \ \ \ &\text{on} \ \partial B_{x_0}(R) \;.
      \end{array}
\label{eq:wrobprob}
\end{equation}
The conclusion for property (C) will follow from the argument preceding Theorem~6.3 of \cite{BH}, 
once we are able to conclude that the sets $\{ w_r < G_R \}$ are a positive distance away from the 
boundary of $B_{x_0}(R)$, for $r \le r_0$ sufficiently small.

First, by comparing with $\max \{ w_r, 0 \}$, note that the minimizer $w_r$ is nonnegative.
It follows from the nondegeneracy theorem, that if we are away from $x_0$ and we are at a free boundary point $y_0$
(i.e. a place where $w_r$ separates from $G_R$), then the difference $G_R - w_r$ obeys the estimate:
\begin{equation}
    \sup_{y \in B_{s}(y_0)} (G_R(y,x_0) - w_r(y)) \geq C_{_{SN}} r^{-n} \min\{s, s^2\} \;,
\label{eq:nondegprop}
\end{equation}
where $C_{_{SN}}$ is the minimum of the two nondegeneracy constants, i.e.
\begin{equation}
     C_{_{SN}} := \min \{ C_{_{ND}}, C_{_{GN}} \} \;.
\label{eq:CMNdef}
\end{equation}
On the other hand, since $G_R(x,x_0)$ is continuous away from $x_0,$ and, for $0 < \delta < r_{_{ND}},$
$\overline{\eta_{2\delta}(\partial B_{x_0}(R))} \cap \overline{B_{x_0}(R)}$ is a compact set, there exists
a positive constant $C_{1}$ such that $0 \leq G_R(x,x_0) \leq C_{1}$ in all of
$\overline{\eta_{2\delta}(\partial B_{x_0}(R))} \cap \overline{B_{x_0}(R)}.$  
Now, we choose $y_0 \in B_{x_0}(R)$ so that $\dist(y_0, \partial B_{x_0}(R)) = \delta$ and 
assume that it belongs to the closure of the positivity set of
$G_R - w_r.$  Applying the nondegeneracy theorems, we see that
\begin{equation}
    C_{1} - \inf_{y \in B_{\delta}(y_0)} w_r(y) \geq \sup_{y \in B_{\delta}(y_0)} (G_R(y,x_0) - w_r(y)) \geq C_{_{SN}} r^{-n} \delta^2\;,
\label{eq:TheGreatnessThatIsC37}
\end{equation}
and this implies:
\begin{equation}
     \inf_{y \in B_{\delta}(y_0)} w_r(y) \leq C_{1} - C_{_{SN}} r^{-n} \delta^2 < 0
\label{eq:HappinessFromC37}
\end{equation}
as soon as $r$ is sufficiently small.
Thus, we have a contradiction if we don't have $G_R(y,x_0) - w_r(y) \equiv 0$ for all $y$ that are
a distance $\delta$ from $\partial B_{x_0}(R).$
By unique continuation, we can be sure that $G_R(y,x_0) - w_r(y) \equiv 0$ on all of
$\overline{\eta_{\delta}(\partial B_{x_0}(R))} \cap \overline{B_{x_0}(R)}.$
Now we are done with the proof of (C) however, as the proof is
the same as the proof starting right after Equation 6.8 in \cite{BH}.

The proofs of (A) and (B) are easier.  The proof of (B) follows from the fact (stated above as IV) that
$\lim_{x \rightarrow x_0} G_R(x,x_0) = + \infty.$  For the proof of (A), the fact that 
$D_{x_0}(r) \subset D_{x_0}(s)$ follows from the claim:
\begin{equation}
    w_s \leq w_r \;.
\end{equation}
Suppose there is a $\tilde{x}$ where $w_s(\tilde{x}) > w_r(\tilde{x})$ and define
$v(x) := w_r(x) - w_s(x).$ By the proof (C) above, we know that $w_s = w_r = G_R(\cdot, x_0)$
in some neighborhood of the boundary $\partial B_{x_0}(R)$. 
Thus, $v$ is identically zero outside of a compact set, and has
a negative minimum at some point $\bar{x}.$  It follows that $w_r(\bar{x}) < G(\bar{x},x_0).$  
Thus, in a neighborhood
of $\bar{x}$ we have
\begin{equation}
    \Delta_g v = \Delta_g w_r - \Delta_g w_s = -r^{-n} - \Delta_g w_s \leq s^{-n} - r^{-n} < 0
\end{equation}
contradicting the maximum principle.  Of course, it is still possible that there is a point
$y_0 \in \partial D_{x_0}(r) \cap \partial D_{x_0}(s)$ which would allow
$$D_{x_0}(r) \subset D_{x_0}(s) \ \ \text{and} \ \ \qhclosure{D_{x_0}(r)} \; \nsubseteq D_{x_0}(s) \;.$$

So, we suppose that there exists a $y_0 \in \partial D_{x_0}(r) \cap \partial D_{x_0}(s).$  By zooming in toward
$y_0$ if necessary, and using the original regularity theorems for the free boundary due to Caffarelli, see for
example \cite[Corollary 3 and Theorem 3]{C0}, we know that there exists a geodesic ball
$B_{z_0}(\rho)$ which satisfies:
\begin{itemize}
    \item[1)] $B_{z_0}(\rho) \subset D_{x_0}(r),$ and
    \item[2)] $\partial B_{z_0}(\rho) \cap \partial D_{x_0}(r) \cap \partial D_{x_0}(s) = y_0.$
\end{itemize}
Considering the function $v(x) := w_r(x) - w_s(x)$ again.  We see that
$v \geq 0$ in $B_{z_0}(\rho),$ and $v(y_0) = 0.$
Since $\Delta_g v = s^{-n} - r^{-n} < 0$ in $B_{z_0}(\rho),$ we can apply the Hopf lemma
to see that $\nabla v(y_0) \ne 0.$  (See for instance Taylor \cite[Chapter 5, Proposition 2.2]{MET}.)
On the other hand, since $y_0$ is in the free boundary for both $w_r$ and $w_s$ we have
\begin{equation}
    \nabla w_r(y_0) = \nabla G(x_0, y_0) = \nabla w_s(y_0) \,,
\label{eq:onthefbs}
\end{equation}
which forces $\nabla v(y_0) = 0,$ giving us a contradiction.
\hfill $\square$

Because of the construction given in the proof above, it makes sense to introduce the notation
$D_{x_0}(r; \mathcal{N})$ to denote the set $\left\{ x \, \big| \, w_r(x) < G_{\mathcal{N}}(x_0,x) \right\}.$
As this notation implies, it may be possible for the {\it noncontact} sets
$\{w_r < G_R \}$ to depend upon the submanifold on which we construct the Green's
function and minimizer $w_r$. However, in the following result, we demonstrate
that the mean value sets are invariant under changes to the submanifold $\mathcal{N}$
if they avoid the submanifold boundary $\partial \mathcal{N}$.   

\begin{lem}[Uniqueness of the Family of Mean Value Sets]    \label{UniqueMethod}
Consider a point $x_0 \subset \subset \mathcal{M},$ $r > 0$ fixed, and $n$--submanifolds 
$\mathcal{N}_1, \mathcal{N}_2$ with boundary, each of which 
contain an open set about $x_0$ and are each compactly contained in $\mathcal{M}.$
We call their respective Green's functions $G_1$ and $G_2,$ and we let
$W_i$ solve $\Delta_g W_i =-r^{-n} \chi_{\{W_i <G_i\}}$ for $i=1,2$ with $W_i \equiv 0$ on
$\partial \mathcal{N}_i.$  Defining the sets 
\begin{equation}
    D_{x_0}^i(r) := D_{x_0}(r;\mathcal{N}_i) = \{ W_i < G_i \} \qquad \text{for $i = 1,2.$}
\label{eq:Drix0def}
\end{equation}
If 
\begin{equation}
     D_{x_0}^1(r) \subset \subset \mathcal{N}_1 \ \ \text{and} \ \ 
     D_{x_0}^2(r) \subset \subset \mathcal{N}_2 \, ,
\label{eq:goodinc}
\end{equation}
then $D_{x_0}^1(r) = D_{x_0}^2(r) \subset \subset \mathcal{N}_1 \cap \mathcal{N}_2.$
\end{lem}

\begin{rmk}
Following from this invariance, we modify our notation for mean value sets whenever we 
can ensure that they avoid the boundary of large enough submanifolds $\mathcal{N}$.
In particular, we will freely use the notation $D_{x_0}(r)$ whenever we assume one can
find a submanifold with boundary on which $\{ w_r < G_{\mathcal{N}}\}$ is uniformly 
bounded away from $\partial \mathcal{N}$. Meanwhile, the notation $D_{x_0}(r;\mathcal{N})$
indicates a construction where the noncontact set may approach the boundary of $\mathcal{N}.$ 
\end{rmk}

\noindent 
\textit{Proof of Lemma\refthm{UniqueMethod}\!\!.}
Choose a compact $\closure{\mathcal{N}}$ with
$\mathcal{N}_1 \cup \mathcal{N}_2 \subset \subset \; \closure{\mathcal{N}}$
and denote its Green's function by $\closure{G}.$ 
For $i =1,2,$ we have that $w=W_i$ minimizes
	$$\int_{\mathcal{N}_i} |\nabla_g w|^2 - 2r^{-n}w$$
among $w \leq G_i$ which are equal to $0$ on $\partial \mathcal{N}_i.$
This implies that $\Delta_g W_i = -r^{-n}\chi_{\{W_i <G_i\}}$ in $\mathcal{N}_i$ and $W_i=0$
on $\partial \mathcal{N}_i.$ Similar equations hold for $\closure{W}$ and $\closure{G}$ where $\closure{G}$ is the
Green's function for $\Delta_g$ on $\closure{\mathcal{N}}$ and $\closure{W}$ will denote the solution of the
corresponding obstacle problem on $\closure{\mathcal{N}}.$

Note that $G_i$ differs from $\closure{G}$ by a harmonic function, specifically, we have
\begin{equation}
      G_i = \; \closure{G} - \; h_i
\end{equation}
where $\Delta_g h_i=0$ and $h_i> 0,$ since $G_i < \; \closure{G}$ by the
maximum principle.  Now we make the following definition:
\begin{equation}
    \tilde{W}_i(x) := \left\{
\begin{array}{rl}
      \closure{G}\!(x,x_0) & \ \text{for} \ x \in \; \closure{\mathcal{N}} \setminus \mathcal{N}_i \\
      W_i(x) + h_i(x) & \ \text{for} \ x \in \mathcal{N}_i \;. \\
\end{array}
           \right.
\label{eq:tildwdef}
\end{equation}
Because $D_{x_0}^i(r) \subset \subset \mathcal{N}_i,$ it is easy to see that
$\tilde{W}_i(x) \equiv \; \closure{G}\!(x,x_0)$ on a neighborhood of $\partial \mathcal{N}_i.$
From there we can verify that $\tilde{W}_i$ satisfies the exact same obstacle problem as
$\closure{W}$ and then, by uniqueness of solutions to the obstacle problem, we get
$\tilde{W}_i \equiv \; \closure{W} \!\!.$  From here we conclude that
	\begin{align*}
		D_{x_0}^i(r) &= \{ W_i < G_i \} \\
                                        &= \{ W_i + h_i < G_i + h_i \} \\
                                        &= \{ \tilde{W}_i < \; \closure{G} \} \\
                                        &= \{ \closure{W} < \; \closure{G} \} \\
                                        &= \; \closure{D}_{x_0} \! (r) \;
	\end{align*}
which proves the result.
\hfill $\square$

Although the previous uniqueness proof breaks down if one of the sets
$D_{x_0}(r, \mathcal{N})$ reaches the boundary of the submanifold $\mathcal{N}$,
the following monotonicity result (with respect to domain inclusion) 
will be useful in this setting. 
\begin{prop}[Noncontact Monotonicity] \label{MonoDProp}
If $\mathcal{N}_1 \subset \mathcal{N}_2 \subset \subset \mathcal{M},$ then for any $x_0 \in \mathcal{N}_1$
and any $r > 0,$ we have the inclusion:
\[
    D_{x_0}(r;\mathcal{N}_1) \subset D_{x_0}(r;\mathcal{N}_2) \;.
\]
\end{prop}

\noindent 
\textit{Proof of Proposition\refthm{MonoDProp}\!\!.}
The idea here is essentially the same as the idea of the proof of the previous lemma.  
Using the same notation as the previous proof, if we let $h$ be the function satisfying
$h := G_2 - G_1$ on $\partial \mathcal{N}_1$ and
$\Delta_g h = 0$ in $\mathcal{N}_1,$ 
then $G_2 \equiv G_1 + h$ in all of $\mathcal{N}_1.$  Furthermore,
letting $W_i$ denote the solution to the obstacle problem:
\begin{equation}
      \begin{array}{rll}
             \Delta_g u &\!\!\!= -r^{-n}\chi_{_{ \{u < G_i \} }} \ \ \ & \text{in} \ \mathcal{N}_i \\
              \ \\
              u &\!\!\!= 0 \ \ \ \ &\text{on} \ \partial \mathcal{N}_i \;.
      \end{array}
\label{eq:wriobprob}
\end{equation}
we see that the functions $w_2 := W_2$ and $w_1 := W_1 + h$ both satisfy
\[
\Delta_g u = -r^{-n}\chi_{_{ \{u < G_2 \} }} \ \ \text{and} \ \ u \leq G_2
\]
within $\mathcal{N}_1.$ While, on $\partial \mathcal{N}_1$ we have $w_2 \leq w_1.$  
The result now follows by the comparison principle for the obstacle problem, 
along with the observation that $\{w_1 < G_2\} = \{W_1 < G_1\}.$  
(See \cite[Theorem 2.7a]{B} for a proof of the comparison principle in 
the Euclidean case that essentialy uses only the weak maximum principle.)
\hfill $\square$

What remains for us to prove now is the existence of a maximal $r_0,$ that allows us to be sure that from a
certain point of view the mean value sets become as large as possible.  Stated differently, we want to
show that we can keep increasing $r$ until the set $D_{x_0}(r)$ either collides with a boundary of
$\mathcal{M}$ or escapes out to infinity.  In the case where our manifolds are compact, we succeed in
proving this conjecture by the time we get to Lemma\refthm{ECS}\!\!.  On the other hand, in the case of
noncompact manifolds we will succeed only partially, and these results are contained in the next section.

\begin{lem}[Convergence of Membranes]  \label{ConMem}
We assume that we are given functions $w_r$ that minimize
\[
	J(w, B_{x_0}(R),r) := \int_{B_{x_0}(R)} |\nabla_g w|^2 - 2r^{-n}w
\]
among functions which vanish on $\partial B_{x_0}(R),$ and are constrained
to be less than or equal to the Green's function with singularity at $x_0.$
Then, as $s \rightarrow r,$
\begin{equation}
    w_s \rightharpoonup w_r \ \text{in} \ W^{1,2}(B_{x_0}(R))
\label{eq:HappyLimit1}
\end{equation}
and
\begin{equation}
    \lim_{s \rightarrow r} ||w_s - w_r||_{C^{\alpha}(\qclosure{B_{x_0}(R)})} = 0
\label{eq:HappyLimit2}
\end{equation}
for some $\alpha > 0.$
Furthermore, as $s \rightarrow \infty$ we have
\begin{equation}
   ||w_s||_{C^{\alpha}(\qclosure{B_{x_0}(R)})} \ \ \text{and} \ \ ||w_s||_{W^{1,2}(B_{x_0}(R))} \rightarrow 0 \;.
\label{eq:atinfty}
\end{equation}
\end{lem}

\noindent 
\textit{Sketch of the proof of Lemma\refthm{ConMem}\!\!.}
 Except for Equation \eqref{eq:atinfty}, the proof is almost
the same as what is given in \cite{AB},
but to make the current work more self-contained, we provide here a sketch.  
First of all, because of the functional
being minimized, it is not hard to show that all of the relevant functions are uniformly 
bounded in $W^{1,2}(B_{x_0}(R)).$
By elliptic regularity, it is then easy to get uniform bounds on the functions in 
$C^{\alpha}.$  So, taking a sequence
$s_n \rightarrow r$ we get a corresponding sequence of functions 
$w_{s_n}$ which converges weakly in $W^{1,2}$
and strongly in $C^{\alpha}$ to a function $\tilde{w}.$  
Now, in order to show that $\tilde{w}$ is in fact $w_r,$ we
use the lower semicontinuity of the Dirichlet integral with respect to $W^{1,2}$ 
convergence along with uniform convergence of $w_{s_n}$ to $\tilde{w}$ to get 
that $J_{r}(\tilde{w}) \leq J_{r}(w_r).$ The claim then follows
by uniqueness of minimizers.

Turning to Equation \eqref{eq:atinfty} we make the following observation:  
If we let $H_r$ be defined to be the solution to the problem
  \begin{equation}
        \begin{array}{rll}
            \Delta_g H_r &\!\!\!= - r^{-n} & \ \text{in}
               \ \ B_{x_0}(R) \vspace{.05in} \\
           H_r &\!\!\!= 0 & \ \text{on} \ \ \partial B_{x_0}(R) \;,
        \end{array}
  \label{eq:Hrdef}
  \end{equation}
then by the maximum principle it is clear that $H_r \geq w_r \geq 0.$  
On the other hand, Schauder theory on
manifolds implies that $H_r \rightarrow 0$ uniformly as $r \rightarrow \infty,$ and then the result
follows from standard estimates on solutions to elliptic PDEs.
\hfill $\square$

We will now show a pair of useful corollaries of this result.  
These corollaries, in some sense, each indicate that
the boundaries of the Mean Value Sets vary continuously.

\begin{coro}[Extension of Mean Value Sets to Greater $r$]   \label{OBANS}
If $D_{x_0}(r) \subset \subset \mathcal{M},$ then there exists an $\epsilon > 0$ 
so that for all $s \in [r, r + \epsilon],$
the set $D_{x_0}(s) \subset \subset \mathcal{M}.$  More precisely, with $w_s$ 
defined as the minimizer of
$J(w, B_{x_0}(R), s)$ the set $\{w_s < G\} \subset \subset \mathcal{M}.$
\end{coro}

\noindent 
\textit{Proof of Corollary\refthm{OBANS}\!\!.}
Without loss of generality we assume that there is a $\delta > 0$ such that the closure of the $2\delta$
neighborhood of $D_{x_0}(r)$ is still compactly contained within $\mathcal{M}.$  We also assume
that $\delta < r_0,$ where $r_0$ is the number from our local nondegeneracy theorem.  (Recall
Theorem\refthm{NonDeg}\!\!.)  Now suppose that there is a sequence $s_k \downarrow r$ such that the
$D_{x_0}(s_k)$ sets are not contained within the $\delta$ neighborhood of $D_{x_0}(r).$  
In this case, there exists a corresponding sequence of points 
$\{y_k\} \subset D_{x_0}(s_k) \cap \partial (\eta_{\delta}(D_{x_0}(r))).$
Since the $y_k$ belongs to $D_{x_0}(s_k),$ we have
\begin{equation}
     w_{s_k}(y_k) < G(y_k, x_0) \;.
\label{eq:ynnotcont}
\end{equation}
Applying Theorem\refthm{NonDeg}to $w_{s_k}$ while observing that $w_{r} \equiv G(\cdot,x_0)$
within $B_{y_k}(\delta)$ for all $k,$ we see that
\begin{equation}
    \sup_{B_{y_k}(\delta)} (w_{r}(y) - w_{s_k}(y)) = 
    \sup_{B_{y_k}(\delta)} (G(y,x_0) - w_{s_k}(y)) \geq
    C_{_{ND}} \delta^{2} s_k^{-n} > 0 \;.
\label{eq:ApplyingND}
\end{equation}
Now, by the uniform convergence of the $w_{s_k}$ to $w_r$, guaranteed by
Equation\refeqn{HappyLimit2}of Lemma\refthm{ConMem}\!\!, we get the desired contradiction.
\hfill $\square$

\noindent
\begin{coro}[Continuous expansion of mean value sets]  \label{CEMVB}
Fix $x_0, y_0 \in \mathcal{M}$ and assume that there exists $0 < s < t$ such that
$y_0 \notin D_{x_0}(s)$ but $y_0 \in D_{x_0}(t).$  Then there exists a unique $R_0 \in (s,t)$
such that $y_0 \in \partial D_{x_0}(R_0).$
\end{coro}

The following proof is essentially identical to the proof in Aryal and Blank \cite{AB} 
which in turn borrows a key idea from the counter-example within Blank and Teka \cite{BT}, but
to keep our work more self-contained we provide a sketch of the proof.

\noindent 
\textit{Sketch of the proof of Corollary\refthm{CEMVB}\!\!.}  Let
$S := \{ t \in \R : y_0 \notin D_{x_0}(t) \}$ 
and then set $R_0 := \sup S.$  We proceed by contradiction. If $y_0 \notin \partial D_{x_0}(R_0),$ then we know that
there exists a small $\rho > 0$ such that either $C:=\overline{B_{y_0}(\rho/2)}$ is completely contained
within $D_{x_0}(R_0)=\{w_{R_0} < G \}$ or it is completely contained within
$D_{x_0}(R_0)^{\complement}=\{w_{R_0}=G\}.$ 

First assume that $C \subset D_{x_0}(R_0).$ Defining $\alpha$ to be the minimum of $w_{R_0}$ in $C,$ Lemma \refthm{ConMem} implies that for small $\delta>0,$ we have $\|w_r-w_{R_0}\|_{L^{\infty}(C)} \leq \alpha/2,$ whenever $|r-R_0|<\delta.$ It follows from the triangle inequality that 
	$$\|w_r - G\|_{L^{\infty}(C)} \geq \|w_{R_0} - G\|_{L^{\infty}(C)}-\|w_r-w_{R_0}\|_{L^\infty(C)}
                                                           \geq \frac{\alpha}{2} >0.$$
Taking $r=R_1:=R_0+\delta/2>R_0,$ the previous inequality tells us that $w_{R_1}< G$ on $C$
implying that $y_0 \in D_{x_0}(R_1),$ contradicting that $R_0=\sup S.$

Now assume that $C \subset D_{x_0}(R_0)^{\complement}=\{w_{R_0}=G\}.$
By nondegeneracy, we know that whenever $r>R_0,$ there exists $z \in C$ and a
constant $\beta >0$ such that $G - w_r(z)> \beta.$ On the other hand, 
Lemma\refthm{ConMem}implies that there exists $\delta >0$ such that 
$\|w_r - w_{R_0}\|_{L^{\infty}(C)} \leq \beta/2.$ Since $w_{R_0} \equiv G$ on $C,$ we have that 
	$$0<\beta \leq \|G - w_r\|_{L^{\infty}(C)}=\|w_r - w_{R_0}\|_{L^{\infty}(C)} \leq \frac{\beta}{2},$$
a contradiction.
\hfill $\square$

\begin{lem}[Exiting compact sets]  \label{ECS}
Fix $R > 0,$ and $x_0 \in \mathcal{M}.$
Then there exists an $r_0 > 0$ such that $\partial D_{x_0}(r_0) \cap \partial B_{x_0}(R)$
is nonempty.
\end{lem}

\noindent 
\textit{Proof of Lemma\refthm{ECS}\!\!.}
This lemma will follow from the claim that as $r$ increases, the minimum distance from
$\partial D_{x_0}(r)$ to $\partial B_{x_0}(R)$ must become arbitrarily small.
Suppose not.  Then there exists a $\delta > 0$ so that within
$S_{\delta} := \eta_{\delta}(\partial B_{x_0}(R)) \cap B_{x_0}(R)$
(i.e. the one-sided $\delta$ neighborhood of $\partial B_{x_0}(R)$) 
we have $w_r \equiv G$ for all $r.$  On $B_{x_0}(R) \setminus S_{\delta}$
we know that $G$ achieves a minimum, $\gamma_{\delta} > 0.$  By the maximum
principle, since $w_r = G$ on $\partial S_{\delta},$ we have
\begin{equation}
   w_r \geq \gamma_{\delta} >0
\label{eq:happycontrad}
\end{equation}
in all of $B_{x_0}(R) \setminus S_{\delta}$ for all $r.$  Now by applying
Equation \eqref{eq:atinfty} we have a contradiction.
\hfill $\square$

\noindent
So with this last lemma, except for showing that $r_0 < \infty,$ 
we can already consider Theorem\refthm{RiemMVT}to be
proven in the case of compact manifolds with boundary, and the $r_0$
from the lemma is the $r_0$ for the theorem. 

\begin{coro}[For Compact Manifolds $r_0 < \infty$]  \label{FCMROLTI}
Let $\mathcal{M}$ be a compact Riemannian manifold with boundary, and let
$x_0 \in \mathcal{M}.$  Then there exists a maximal $r_0 < \infty$ such that
$D_{x_0}(r_0) \subset \mathcal{M}$ and
$\partial D_{x_0}(r_0) \cap \partial \mathcal{M}$ is nonempty.
\end{coro}

\noindent 
\textit{Proof of Corollary\refthm{FCMROLTI}\!\!.}
Suppose not.  Then the set
\[
     D := \bigcup_{r > 0} D_{x_0}(r)
\]
has a positive distance from $\partial \mathcal{M}$ and so by applying
Lemma\refthm{ECS}to a compact $K$ satisfying
\[
     D \subset K \subset \subset \mathcal{M}
\]
we get a contradiction.
\hfill $\square$

\section{Results for Strongly Nonparabolic Manifolds}  \label{Nonpara}

The results of the previous sections offer a foundation for the local theory of the MVT on complete 
Riemannian manifolds.  On the other hand,
in the case where $\mathcal{M}$ is not compact some obvious questions remain.  The most pressing and
immediate issue seems to be a statement that $r_0$ ``gets as large as possible.''  Indeed, one could
hope for a statement that the mean value balls can grow until either they hit a finite boundary or until
they escape to infinity.  Following Lei Ni we make the following definition \cite[Definition 2.2]{N}:
\begin{defi}[Strongly Nonparabolic Riemannian Manifolds]  \label{SNParab}
We call the unbounded Riemannian manifold $\mathcal{M}$ strongly nonparabolic if there is a minimal
positive Green's function $G_{\infty}(x,y)$ which satisfies:
\begin{equation}
     \lim_{y \rightarrow \infty} G_{\infty}(x,y) = 0
\label{eq:GgtZ}
\end{equation}
for fixed $x \in \mathcal{M}.$
\end{defi}
Ni discusses this definition and shows a variety of classes of manifolds which satisfy this property \cite{N}.
Indeed, in addition to the trivial observation that $\R^n$ works for $n \geq 3,$ Ni mentions the following:
\begin{itemize}
    \item[a)] Li and Yau \cite{LY} prove that if $\mathcal{M}$ is a manifold with nonnegative Ricci curvature,
then $\mathcal{M}$ is nonparabolic if and only if
\begin{equation}
   \int_{r}^{\infty} \frac{\tau}{\text{Vol}(B_{x_0}(\tau))} d\tau < \infty \;.
\label{eq:eqtsnp}
\end{equation}
There is also an estimate given which shows that for this case nonparabolic implies strongly nonparabolic.
    \item[b)] If $\mathcal{M}$ satisfies the Sobolev Inequality:
\begin{equation}
           \left( \int_{\mathcal{M}} f^{\frac{2\nu}{\nu -2}} \, d\mu \right)^{\frac{\nu - 2}{2\nu}}
    \leq A\int_{\mathcal{M}} |\nabla f|^2 d\mu
\label{eq:HappySOB}
\end{equation}
for some $\nu > 2, A > 0,$ and all smooth $f$ which are compactly supported, then $\mathcal{M}$ is
strongly nonparabolic.
    \item[c)] Finally Li and Wang \cite{LW} give examples related to the spectrum of the Laplace
operator \cite{LW}. Further, Ni shows that if $\mathcal{M}$ has a positive lower bound on the spectrum 
of the Laplace operator, if $Ric \geq -(n - 1)g,$ and finally if the volume of the unit ball around 
every point is bounded uniformly from below by a positive constant, then $\mathcal{M}$ is strongly nonparabolic
\cite[Proposition 2.5]{N}.
\end{itemize}

In any case when $\mathcal{M}$ is strongly nonparabolic we have the following boundedness result:
\begin{theo}[Mean value sets are bounded]   \label{MVSAB}
Assume that $\mathcal{M}$ is complete, unbounded without boundary, and strongly nonparabolic, and
assume further that $r_{_{ND}}$ and $C_{_{ND}}$ are bounded from below at every point of $\mathcal{M}.$
Then given any $x_0 \in \mathcal{M}$ and $r > 0,$ there exists an $s > 0$ such that
\[
   D_{x_0}(r;B_{x_0}(s)) \subset \subset B_{x_0}(s) \;.
\]
Of course in this case we can say that $D_{x_0}(r;B_{x_0}(s)) = D_{x_0}(r),$ and as a consequence
$r_0 = \infty.$
\end{theo}
\noindent 
\textit{Proof of Theorem\refthm{MVSAB}\!\!.}
The proof will use the same key idea as the proof we gave earlier for Lemma\refthm{RiemMVTagain}\!\!.
As in that proof we will abuse notation by using ``$B_{x_0}(s)$'' to mean a smooth approximation to
$B_{x_0}(s)$ when necessary.  Also, since we will always be evaluating our Green's functions at $x_0$
as one of the points, we will suppress that within our notation below.

Since $\mathcal{M}$ is strongly nonparabolic, we will let $G_{\infty}$ denote the minimal Green's function
satisfying Equation \eqref{eq:GgtZ}.  Given any $\delta > 0,$ there exists an $s_0$
such that $G_{\infty}(x) \leq \delta$ in all of $B_{x_0}(s_0)^{\complement}.$
Define $s_1 := s_0 + 2r_{_{ND}},$ and let $G_{1}(x)$ denote the
Green's function on $B_{x_0}(s_1)$ we note the following:
\begin{itemize}
    \item[a)] $0 \leq G_{1}(x) \leq G_{\infty}(x)$ in $B_{x_0}(s_1),$
    \item[b)] $G_{\infty}(x) \leq G_{1}(x) + \delta$ in $B_{x_0}(s_1) \setminus B_{x_0}(s_0),$ and
    \item[c)] $h(x) := G_{\infty}(x) - G_{1}(x)$ is a harmonic function in $B_{x_0}(s_1)$
which satisfies the estimate:
\begin{equation}
    ||h||_{L^{\infty}(B_{x_0}(s_j))} \leq \delta
\label{eq:smallnessofthehos}
\end{equation}
by the weak maximum principle.
\end{itemize}

Now let $w_{1}$ denote the minimizer of
\begin{equation}
       J(w, B_{x_0}(s_1), r) := \int_{B_{x_0}(s_1)} |\nabla_g w|^2 - 2r^{-n} w \;,
\label{eq:fctldefnow}
\end{equation}
among functions in $W^{1,2}_0(B_{x_0}(s_1))$ which are less than or equal to
the Green's function, $G_{1}(x).$  As before, this minimizer solves the obstacle problem:
\begin{equation}
      \begin{array}{rll}
             \Delta_g u &\!\!\!= -r^{-n}\chi_{_{ \{u < G_1 \} }} \ \ \ & \text{in} \ B_{x_0}(s_1) \\
              \ \\
              u &\!\!\!= 0 \ \ \ \ &\text{on} \ \partial B_{x_0}(s_1) \;.
      \end{array}
\label{eq:wrobprobnow}
\end{equation}
Now suppose that $y_0 \in \partial B_{x_0}(s_0 + r_{_{ND}})$ and $w_{1}(y_0) < G_{1}(y_0).$
Then we can apply Theorem\refthm{NonDeg}to the function $G_1 - w_1$ at $y_0$ in order to state:
\begin{equation}
    \sup_{B_{y_0}(r_{_{ND}})} (G_1(x) - w_1(x)) \geq r^{-n} C_{_{ND}} (r_{_{ND}})^2
\label{eq:estimateonsup}
\end{equation}
On the other hand,
\begin{equation}
    \sup_{B_{y_0}(r_{_{ND}})} (G_1(x) - w_1(x)) \leq \delta,
\label{eq:otherestimateonsup}
\end{equation}
and so if $\delta$ is sufficiently small, we get a contradiction.  Therefore, we must have
$w_1 \equiv G_1$ in $B_{x_0}(s_0 + r_{_{ND}})^{\complement}.$  Defining $\delta$ to be
such a number, we can take $s$ for the theorem to be the constant $r_{_{ND}}$ plus the
$s_0$ determined by the requirement: $G_{\infty}(x) \leq \delta$ in all of $B_{x_0}(s_0)^{\complement}.$
\hfill $\square$
\begin{rmk}[Weaker Condition]  \label{Weaker}
It follows from the proof that there is a weaker condition under which
Theorem\refthm{MVSAB}holds.  Indeed, if we let $r_{_{ND}}(x)$ and
$C_{_{ND}}(x)$ denote the values of $r_{_{ND}}$ and $C_{_{ND}}$ at the point $x \in \mathcal{M},$
then we can replace all of the hypotheses of the theorem with the simple requirement that a Green's
function $G_{\infty}$ exists and
\begin{equation}
    \lim_{y \rightarrow \infty} \frac{G_{\infty}(x,y)}{C_{_{ND}}(y) r_{_{ND}}(y)^2 } = 0
\label{eq:MagicCondition}
\end{equation}
for any $x \in \mathcal{M}.$  Of course, since all of the statements above can be made for $\R^2$
which is obviously not a strongly nonparabolic manifold, and which does not satisfy
Equation \eqref{eq:MagicCondition}, it is clear that there must be weaker
conditions than even this one.
\end{rmk}

\bibliographystyle{plain}
\bibliography{Draft}
\end{document}